\documentclass[11pt]{article}

\usepackage{amssymb,latexsym,amsfonts,verbatim,amscd}

\newtheorem{Def}{Definition}[section]
\newtheorem{Thm}[Def]{Theorem}
\newtheorem{Lem}[Def]{Lemma}
\newtheorem{Prop}[Def]{Proposition}
\newtheorem{Cor}[Def]{Corollary}
\newtheorem{Rem}[Def]{Remark}

\newtheorem{Fac}[Def]{Fact}

\newtheorem{Quest}[Def]{Question}
\newtheorem{Quests}[Def]{Questions}
\def\telos{\hfill$\dashv$}

\begin{document}
\baselineskip=17pt

\title{Consequences of Vop\v{e}nka's Principle over weak set theories}
\author{Athanassios Tzouvaras}

\date{}
\maketitle

\begin{center}
 Department  of Mathematics\\  Aristotle University of Thessaloniki \\
 541 24 Thessaloniki, Greece \\
  e-mail: \verb"tzouvara@math.auth.gr"
\end{center}

$$\mbox{\em  Dedicated to the memory of Petr Vop\v{e}nka}$$

\vskip 0.2in

\abstract{It is shown that  Vop\v{e}nka's Principle (VP) can restore almost the entire ZF over a  weak fragment of it. Namely, if EST is the theory consisting of the axioms of Extensionality,  Empty Set, Pairing, Union, Cartesian Product, $\Delta_0$-Separation and Induction along $\omega$, then ${\rm EST+VP}$ proves the axioms of Infinity, Replacement (thus also Separation) and Powerset. The result was  motivated by previous results in \cite{Tz14}, as well as by H. Friedman's \cite{Fr05}, where a  distinction is made among   various  forms of VP.  As a corollary,   ${\rm EST}+$Foundation$+{\rm VP}$=${\rm ZF+VP}$,  and ${\rm EST}+$Foundation$+{\rm AC+VP}={\rm ZFC+VP}$.  Also it is shown that the Foundation axiom is independent from  ZF--\{Foundation\}+${\rm VP}$. It is open  whether the Axiom of Choice is independent from ${\rm ZF+VP}$. A very weak form of choice follows from VP and some similar other forms of choice are introduced.}

\vskip 0.2in


\vskip 0.2in


\section{Introduction}

Vop\v{e}nka's Principle (henceforth abbreviated VP) is mainly known as a (very) large cardinal axiom (see \cite{Je03}). Also  several other implications of the principle have been  proved long ago, especially in category theory (see \cite{AR94}). Recently there  has been a revived  interest in VP through new set theoretic proofs of category theoretic results (see Bagaria and Brooke-Taylor \cite{BBT13}). Also Brooke-Taylor showed  in \cite{BT11} the  relative consistency of VP with almost all usual ZFC-independent statements, like $GCH$ and  diamond principles (see  the Introduction of \cite{BT11}).

In all of the above results the  underlying theory is  ${\rm ZFC}$.
In  contrast, the aim of the present paper  is to reveal  a still different capability  of VP: the capability to restore the most basic axioms of ZF, namely  Replacement (thus also Separation) and Powerset, as well as Infinity (if the latter is missing),   when added to a suitable weak set theory. By ``weak set theory'' we generally  mean the following.

\begin{Def} \label{D:weaktheory}
{\em A} weak set theory {\em is one that   does not include the axioms of Powerset and  Replacement.}
\end{Def}

Various weak systems of set theory have been considered in the literature. Perhaps the most well-known of them is  Kripke-Platek's theory KP on which the theory of admissible sets is based (see \cite{De84}, p. 48). A system weaker than KP is Devlin's  Basic Set Theory (BS) used in \cite{De84} (see \cite{De84}, p. 36).  An extensive and detailed treatment of an array of weak systems, among them Devlin's BS, can be found in \cite{Mat06}.

A weak set theory may or may not include Infinity. Also it may or may not be a fragment of ZF. For example the theory LZFC (``local ZFC'') of \cite{Tz10} and \cite{Tz14} is weak, proves Infinity,  but is not a fragment of ZFC. On the other hand, the system EST introduced  in section 2 below,  as a ground theory for VP, is a fragment of ZF but does not include Infinity. Yet EST+VP proves this axiom. It is worth pointing out  that EST, even augmented with Infinity, is weaker than BS because of lack of Foundation. {\em A fortiori} it is weaker than KP.

Throughout we shall  refer to the well-known axioms of ZFC with their usual names, and without further explanations. These are:  Extensionality, Empty Set, Pairing, Union, Powerset, Infinity, Separation, Foundation, Replacement and Choice. Sometimes it is convenient to denote them by  abbreviations, especially within theorems. Specifically we often  write $Ext$ for Extensionality, $Pair$ for Pairing, $Pow$ for Powerset, $Sep$ and $\Delta_0$-$Sep$ for Separation and $\Delta_0$-Separation, respectively, $Rep$ and $\Delta_0$-$Rep$ for Replacement and  $\Delta_0$-Replacement, respectively, $Found$ for Foundation, $Inf$ for Infinity and AC for Choice. Another weak axiom that will be used below is Cartesian Product, abbreviated $CartProd$, which says that for any sets $x,y$, $x\times y$ is a set.

The  capability of  VP to restore the axioms of Replacement and Powerset  was first noticed  in \cite{Tz14}. In that paper we showed that if we add  VP to a strengthened  variant of LZFC,  then  Replacement and  Powerset  are recovered. Later on, when we came across \cite{Fr05},  where seven  variants of VP are given in apparently decreasing strength but still equivalent over ZFC, we realized that what was used in \cite{Tz14} was  not the full principle VP but only a weaker form, denoted  ${\rm VP}_4$ in Friedman's list. Moreover  we saw that using VP instead of ${\rm VP}_4$,  the result of \cite{Tz14}  holds    for the theory LZFC itself rather than a  strengthened variant of it. That led us to focus on what VP can prove over a weak fragment of ZF rather than  LZFC.

\vskip 0.1in

The content of the paper is as follows. In  section 2.1 we introduce the variants of VP, especially the general one VP and the weaker one ${\rm VP}_4$, and describe their basic difference. Also we outline the way in which VP acts as a set-existence principle. In section 2.2 we introduce the weak theory EST. In section 2.3 we show that VP is expressible in EST. Section 2.4 contains the main results of section 2.
Namely that  EST+VP proves  Infinity, Replacement and Powerset.

In section 3 we prove that Foundation is independent from ZF+VP, that is, if ZF+VP is consistent, then so is ${\rm ZF}_0+{\rm VP}+\neg Found$, where ${\rm ZF}_0={\rm ZFC}-\{Found\}$.

In section 4  the question whether AC is independent from ZF+VP is raised. The question remains  open. In particular it is open whether the question can be settled with the help of symmetric and permutation models. It is also observed that a very weak form of choice follows from VP. This gives the chance to introduce some  similar other forms of choice whose relative strength over ZF, as well as over ZF+VP, is also open.

\section{Vop\v{e}nka's Principle over some weak set theories}

\subsection{VP and its variants}
Let ${\cal L}=\{\in\}$ be the language of set theory.  Given a formula $\phi(x)$ of ${\cal L}$ in one free variable, let $X_\phi$ denote the extension $\{x:\phi(x)\}$ of $\phi(x)$. As usual we refer to $X_\phi$ as ``classes''.
Vop\v{e}nka's Principle  is a statement that quantifies over classes, so cannot be formulated in ZF as a single axiom; it can be formulated however as an axiom-scheme.  Clearly, for every $\phi$ ``$X_\phi$ is a proper class'' is a first-order sentence. Therefore so are also the  statements:

\vskip 0.1in

 $({\rm VP}_\phi$) \ \ {\em If $X_\phi$ is a proper class  of $L$-structures, for some  \\  \hspace*{0.7in} first-order language $L$,  then there  are distinct ${\cal M}$, ${\cal N}\in X_\phi$
\\ \hspace*{0.7in} such that  ${\cal M}\precsim {\cal N}$,}

\vskip 0.1in

\noindent where  ${\cal M}\precsim {\cal N}$ means that there is an elementary embedding $f:{\cal M}\rightarrow {\cal N}$. Let   $${\rm VP}=\{{\rm VP}_\phi:\phi(x) \ \mbox{a formula of} \ {\cal L}\}.$$
[Notice that we refer to the arbitrary first-order languages mentioned in the statement of ${\rm VP}_\phi$ above by the symbol ``$L$'',  to avoid confusion  with the fixed language ${\cal L}=\{\in\}$ of set theory. Of course ${\cal L}$ is {\em one} of those $L$, but  a very specific one.  This practice is followed throughout the paper.]

Let us refer to  the  above formulation  of  ${\rm VP}_\phi$ as being ``direct''. We have also the   contrapositive formulation: ``If $X_\phi$ is a  class  of $L$-structures, for some  first-order language $L$,  and  there are no  distinct ${\cal M},{\cal N}\in X_\phi$ such that ${\cal M}\precsim {\cal N}$, then  $X_\phi$ is not a proper class.'' But ``$X_\phi$ is not a proper class'' means exactly that $X_\phi$ is a set.
This latter formulation  enables ${\rm VP}_\phi$ to act as a set-existence principle: It says that   ``if  such and such is the case about $\phi$ and $X_\phi$, then  $X_\phi$ is a {\em set}''.

The direct formulation of VP is suitable for ZF and ZFC, where it is usually easy to decide  whether its  premise  is satisfied, that is, whether a class $X_\phi$ of structures is proper. In contrast,  the contrapositive formulation  makes  VP suitable for weak set theories. The reason is that such theories have poorly defined universes,  where it is often  unclear which classes $X_\phi$ are sets and which are proper ones. So in weak set theories we are more in need of principles entailing that such and such  classes $X_\phi$ are sets. For example, in a weak set theory,  where Replacement and Powerset are missing, such  instances  of VP might  be used to establish  that the classes $\{x:(\exists y\in A)(F_\psi(y)=x)\}$, where $F_\psi(y)=x$ is a functional relation, or $\{x:x\subseteq A\}$,  for any set $A$, are actually sets, thus proving the aforementioned axioms. As already said above, we first noticed this fact   when working on \cite{Tz14}. Motivated by results in that paper and also by Friedman's \cite{Fr05},   we pursued the above  idea more systematically and  showed  that in fact,  by the  contrapositive action  of VP,  the three pillar axioms of ZF,  Powerset, Replacement and Infinity, can be established over a very weak fragment of ZF.

H. Friedman gave in \cite{Fr05} seven variants of VP, namely VP above plus six  weaker forms ${\rm VP}_i$, $1\leq i\leq 6$, almost in decreasing strength, each of them resulting either from  a narrowing of the range of first-order languages that can be engaged in the scheme,  or from  replacing ``elementary embedding'' with ``embedding'', or from both. Specifically, ${\rm VP}_1$ results from VP by restricting the term  ``first-order language''  to  ``language of finite relational type''. ${\rm VP}_2$  results from ${\rm VP}_1$ by replacing ``elementary embedding'' by ``embedding''. ${\rm VP}_3$ results from ${\rm VP}_2$ by replacing ``language of finite relational type'' by ``language with a single binary operation''. ${\rm VP}_4$ results from ${\rm VP}_3$ by replacing ``language with a single binary operation'' with ``language with a single binary relation''. ${\rm VP}_4$ is essentially the variant used in \cite{Tz14}, so we omit the other variants of VP. More precisely, let

\vskip 0.1in

$({\rm VP}_{4,\phi}$) \ \ {\em If $X_\phi$ is a proper class  of $L$-structures of  a language $L$
\\  \hspace*{0.7in}
with a single binary relation ${\bf R}$,  then there  are distinct
\\ \hspace*{0.7in} ${\cal M},{\cal N}\in X_\phi$ such that ${\cal M}$ is embeddable into ${\cal N}$.}

\vskip 0.1in

\noindent Let also
$${\rm VP}_4=\{{\rm VP}_{4,\phi}: \phi(x) \ \mbox{is a formula of  ${\cal L}$}\}.$$
As  Friedman states in  \cite{Fr05},  all these seven variants  are equivalent over the theory of classes Neumann-G\"{o}del-Bernays with Choice (NGBC), when ${\rm VP}_i$ are formulated in the language of the theory of classes. That means that  the above formulations of ${\rm VP}_i$ in ZF,  as schemes, are also  equivalent over ZFC. The equivalence of ${\rm VP}_i$ should have been well-known rather  since the start, at least  among category theorists. Also  as  indicated in the Historical Remarks on pp. 278-279 of \cite{AR94}, the original formulation of Vop\v{e}nka's Principle was ${\rm VP}_4$ rather than VP.  However these variants need  {\em not} be   equivalent over weaker theories, like the one we deal with  below (see Corollary \ref{C:Pultr} for the difference between VP and ${\rm VP}_4$ over EST, and the comments following that).

Although in this paper we have made the shift from the weak variant ${\rm VP}_4$ to the most general one VP, still in the proofs of our  main results below we do not need  the full generality of VP. To be precise: the main difference between ${\rm VP}_4$ and VP, as used below,  is that the first-order languages  involved in VP  may contain, for our purposes, in contrast to those  involved in ${\rm VP}_4$, an {\em arbitrary  set of constants.} On the other hand, for our purposes, the languages in VP need not contain more than one binary and one unary predicate. Moreover, for the derivation of Replacement and  Powerset it suffices that  VP talks just about {\em embeddings} rather than elementary embeddings. (See Remark \ref{R:reduce} below.)

Having made the distinction between VP and ${\rm VP}_4$, we can now describe the  result of \cite{Tz14} that has largely motivated this paper. (Although in \cite{Tz14} we refer to VP, what we actually use is ${\rm VP}_4$.) To formulate it we need first an  old theorem of  P. Vop\v{e}nka, A. Pultr and Z. Hedrl\'{i}n \cite{VPH65}, that we abbreviate V-P-H, the proof of which was crucial for the result of \cite{Tz14}.

Given a set $A$ with a binary relation $R\subseteq A\times A$, we refer to the ordered pair  $\langle A,R\rangle$ as a {\em graph}. Given a graph $\langle A,R\rangle$, a mapping $f:A\rightarrow A$ is  an {\em endomorphism} if for all $x,y\in A$, $\langle x,y\rangle\in R$ implies $\langle f(x),f(y)\rangle\in R$. The graph $\langle A,R\rangle$ is said to be {\em rigid}  if the only endomorphism of $\langle A,R\rangle$ is the identity.

The V-P-H theorem is  the following:

\begin{Thm} \label{T:keyVopenka}
{\rm (V-P-H \cite{VPH65})} {\rm (ZFC)}
For every infinite set  $A$, there is a binary  relation $R\subset A\times A$ such that the graph $\langle A,R\rangle$ is rigid.
\end{Thm}

Obviously  the preceding theorem refers to structures with a single binary relation only. Also rigidity is a property that refers to (lack of even)  {\em endomorphisms} rather than  embeddings (let alone elementary embeddings). This is why what is proved in   \cite{Tz14} is in essence the following.

\begin{Thm} \label{T:main}
{\rm (\cite[Theorem 6.3]{Tz14})} If   $T$ is a theory such that ${\rm LZFC}\subseteq T$ and $T\vdash $V-P-H, then $T+{\rm VP}_4$ proves Replacement and Powerset, that is, ${\rm ZFC}\subseteq T+{\rm VP}_4$.
\end{Thm}
[This  formulation differs from that of \cite{Tz14} only in that we mention ${\rm VP}_4$ in place of  VP.\footnote{The precise definition of the theory LZFC is not needed here. It suffices to say that it is a weak set theory in the sense of definition \ref{D:weaktheory}. It is   not a fragment of ZF, but it satisfies Infinity (in the sense that $\omega$ is the least inductive set), AC, Cartesian Product, and others. For later use we note also that LZFC is much stronger than the theory EST introduced below.}]

In Theorem \ref{T:main}, the assumption that the theory $T$ proves  V-P-H guarantees the existence  of  a rigid binary relation $R$ on every set $A$. Then  applying ${\rm VP}_4$ to structures suitably equipped with such an $R$, it is shown that  ${\cal P}(A)$ and $F``A$ are not proper classes. The  variant ${\rm VP}_4$ suffices for this purpose since  we need to employ a language with a  binary relation $\textbf{R}$ only, as well as just {\em endomorphisms} instead of embeddings or  elementary embeddings.

The basic observation that led from Theorem \ref{T:main} above to the results of the present article is that, in the absence of V-P-H, the rigidity property can be alternatively guaranteed by employing an  infinity of constants in the language of the structure, namely a constant $\textbf{c}_a$ for each  element  $a\in A$. Such languages are  allowed in VP though not in ${\rm VP}_4$.

\subsection{A weak fragment of ZF for expressing VP}
Below we define the  weak  fragment of ZF called EST (for Elementary Set Theory). First let us set
$${\rm EST}_0=\{\mbox{\em{Ext,  Emptyset, Pair, Union, CartProd, $\Delta_0$-Sep}}\}.$$
For any sets $x,y$, the ordered pair $\langle x,y\rangle$ is defined in ${\rm EST}_0$ as usual, that is, as the set  $\{\{x\},\{x,y\}\}$. For any given sets $M,N$, $M\times N=\{\langle x,y\rangle:x\in M, y\in N\}$ is a set by $CartProd$. Binary relations between $M$ and $N$ are defined as sets $R\subseteq M\times N$.   Functions $f:M\rightarrow N$ are defined  as special binary relations $f\subseteq M\times N$. Throughout the symbol $f$  ranges over functions, so $(\exists f)(\cdots)$ abbreviates $(\exists f)(f \ \mbox{is a function} \wedge \cdots)$.

Next we define  the class  $\omega$ of natural numbers, which also need not be a set. Let $Tr(x)$ denote the  predicate ``$x$ is a transitive set''. Let us also define  the predicates
$$Ord(x):=Tr(x) \wedge (\forall y,z\in x)(y\in z\vee y=z\vee z\in y),$$
$$Succ(x):=(\exists y\in x)(x=y\cup\{y\}),$$
$$Nat(x):=Ord(x) \wedge (\forall y)[(y\in x\vee y=x)\rightarrow (y=\emptyset\vee Succ(y))].$$
Finally let us set
$$\omega=\{x:Nat(x)\}.$$
We call $\omega$ the {\em class of natural numbers.}
Throughout, in writing formulas of ${\cal L}$, it is convenient to use the notation   $x\in\omega$, although $\omega$ is a class in general, as an abbreviation of the predicate $Nat(x)$. In particular, $(\exists x\in\omega)\phi$ and $(\forall x\in\omega)\phi$ stand for $(\exists x)(Nat(x)\wedge \phi)$ and $(\forall x)(Nat(x)\rightarrow \phi)$, respectively. Also as usual the letters $m$ and $n$ will range over elements of $\omega$.

Further we need {\em induction} to hold along $\omega$, that is, that every non-empty subclass $X\subseteq \omega$  has a least element. For that purpose  an additional axiom is needed. This is the Induction scheme, $Ind(\omega)$, given below, which is almost identical to the Induction scheme of Peano Arithmetic. Namely, for every formula $\phi(x)$ of ${\cal L}$ such that $X_\phi\subseteq\omega$ (that is, $(\forall x)(\phi(x)\rightarrow Nat(x))$), let $Ind_\phi(\omega)$ denote the formula

$$(Ind_\phi(\omega))
\hspace{.5\columnwidth minus .5\columnwidth}
[\phi(\emptyset) \wedge (\forall x\in \omega)(\phi(x)\rightarrow \phi(x\cup\{x\})]\rightarrow (\forall x\in \omega)\phi(x) .\hspace{.5\columnwidth minus .5\columnwidth} \llap{}$$

Let also
$$Ind(\omega)=\{Ind_\phi(\omega):X_\phi\subseteq \omega\}.$$
[Equivalently $Ind(\omega)$ says that every nonempty $X_\phi\subseteq \omega$ has a $\in$-least element.]
Finally we set
$${\rm EST}={\rm EST}_0+Ind(\omega).$$
EST is the weak theory that will be used below as a base theory for VP.

\begin{Rem} \label{R:Ind}
{\em The definition of $\omega$ above is as in \cite[pp. 468f]{BM77} (see also \cite{BF93}), except that in \cite{BM77} the predicate $Ord(x)$ says that    ``$x$ is well ordered with respect to $\in$'', while in our definition of $Ord(x)$,  ``$x$ is ``linearly  ordered with respect to $\in$''. Note the following:

(i) A difference between the two properties is that ``linearly ordered'' is $\Delta_0$ while ``well ordered'' is not (it is $\Pi_1$). This fact  will be needed in the proof of Theorem \ref{T:inf} below.  On the other hand,  with the help of $Ind(\omega)$, one can easily prove that for every $x\in \omega$, $x$ is indeed well ordered with respect to $\in$.

(ii) If Foundation were available, the properties ``$x$ is  linearly ordered'' and ``$x$ is  well ordered'' (with respect to $\in$), would be equivalent.

(iii) However even if Foundation were available, ${\rm EST}_0+Found$ could not prove $Ind(\omega)$. That would need in addition Separation or Replacement. That is, ${\rm EST}_0+Found+Sep\vdash Ind(\omega)$.}
\end{Rem}

With the help of $Ind_\phi(\omega)$ one can prove all basic facts about natural numbers. Some of them needed below are the following:

\begin{Fac} \label{F:Natfacts}

(i) For any $x,y\in \omega$, if $x\in y$ then $x\subsetneq y$ and $x$ is an initial segment of $y$.

(ii) If $x,y\in \omega$ and $x\subsetneq y$, then there is no injection $f:y\rightarrow x$.

(iii) For any $x,y\in \omega$, $x\in y\vee y\in x \vee x=y$.

(iv) If $x\in \omega$ and $x=y\cup\{y\}$, $y$ is the greatest element of $x$.

\end{Fac}

Having defined ordered pairs and natural numbers, ordered $n$-tuples, for $n\in \omega$, can be defined as usual by induction. Namely, for every $n>2$, $\langle x_0,\ldots,x_{n+1}\rangle=\langle\langle x_0,\ldots,x_n\rangle,x_{n+1}\rangle$. Moreover, by the axiom $CartProd$, for  $n\in\omega$, $n>0$, and every $n$-tuple of sets $M_0,\ldots,M_{n-1}$,  the set $M_0\times\cdots\times  M_{n-1}=\{\langle x_0,\ldots,x_{n-1}\rangle:x_i\in M_i\}$ is a set.

Below we shall need also the set of functions $f:n\rightarrow M$,  denoted ${}^nM$, for any set $M$ and any $n\in\omega$. Since Replacement is not available in EST we shall define the elements of ${}^nM$  in a slightly different way. Namely, while an $f\in {}^nM$ has typically the form $$f=\{\langle 0,x_0\rangle,\ldots,\langle n-1,x_{n-1}\rangle\},$$ obviously we can identify the latter set with the $n$-tuple
$$\langle \langle 0,x_0\rangle,\ldots, \langle n-1,x_{n-1}\rangle\rangle,$$
which is an element of the Cartesian product $$(\{0\}\times M)\times\cdots \times (\{n-1\}\times M).$$
The latter is a set as we saw above, so we can define
$${}^nM:=(\{0\}\times M)\times\cdots \times (\{n-1\}\times M).$$
As usual we shall let letters $g,h$ range over elements of ${}^nM$ and we shall write $h(i)=x$ to denote the fact that $\langle i,x\rangle$ is the $i$-th  element of the $n$-tuple $h$.

As already mentioned in the Introduction,  EST, even augmented with Infinity, is weaker than Devlin's system BS. In connection with the  definition of the $Sat$ predicate in EST that will be given in the next section,  the referee kindly informed me that Devlin's  BS, though stronger than EST, is famously  insufficient for the purpose he introduced it, including defining the  satisfaction predicate.\footnote{A discussion on this issue can be found at

http:// mathoverflow.net/questions/77734/devlins-constructibility-as-a-resource.}
(The flaws of BS with respect to this point are discussed and remedied in \cite[\S 10]{Mat06}.)  The main reason that EST succeeds where BS fails  appears to be the unorthodox definition of the set ${}^nM$ given above. Unexpectedly enough, this unusual yet legitimate formalization of ${}^nM$ is all we need to make  things work.

\subsection{Expressibility of VP in EST}
From now on we work in ${\rm EST}$, with language ${\cal L}=\{\in\}$.
In this subsection we  show  how the concepts required  for the formulation of VP can be defined in ${\rm EST}$.
The formulation of VP (at least for the needs of present article) requires the following:

(a) The definition of a language $L_A$, for every set $A$, that contains at most one  unary and one binary relation symbol, but contains a constant  $\textbf{c}_a$  for each $a\in A$ and infinitely many variables. Due to the lack of Replacement  $L_A$ need not be a set.

(b) The definition of terms,   formulas and sentences of $L_A$,  in a way that these classes of objects are {\em inductive,} so that one can  prove inductively the usual facts about these syntactic objects. Again the classes of terms, formulas and sentences  need not be sets.

(c) The definition of $L_A$-structures for every language $L_A$.

(d) The definition of satisfaction relation ``${\cal M}\models \sigma(x_0,\ldots,x_{n-1})$'', for any $L_A$-structure ${\cal M}$,  any formula $\sigma(v_0,\ldots,v_{n-1})$ of $L_A$, with free variables $v_0,\ldots,v_{n-1}$, and any $n$-tuple of elements $x_0,\ldots,x_{n-1}$ of $M$. In view of this, the relation of elementary embeddability ${\cal M}\precsim {\cal N}$ (or simple embeddability) between $L_A$-structures is immediately defined.

We show below how  definitions (a)--(d) can be implemented in EST. Since the definability of the satisfaction relation ${\cal M}\models \sigma(x_0,\ldots,x_{n-1})$ is crucial for the expressibility of VP in EST, we shall give explicitly the necessary definitions below.

\begin{Def} \label{D:prima}
{\em For every $n\in \omega$, let $v_n:=\langle 0,n\rangle$. $v_n$ is the $n$-th} variable {\em (of every language). For every set $a$, let $\textbf{c}_a:=\langle 1,a\rangle$. The sets $\textbf{c}_a$ are called} constants. {\em Also let us identify the (sufficient) logical symbols $\equiv$, $\neg$, $\wedge$ and $\exists$ with  elements of $\omega$ as follows:
$\equiv:=2$, $\neg:=3$, $\wedge:=4$, $\exists:=5$. Finally let
$\textbf{U}:=6$ and $\textbf{R}:=7$. The symbols  $\equiv$, $\textbf{U}$ and $\textbf{R}$ are referred to as} predicates. {\em $\equiv$ and $\textbf{R}$ are} binary {\em predicates, while $\textbf{U}$ is} a unary  {\em one.  For every set $A$, let $$L_A=\{\equiv,\neg,\wedge,\exists\}\cup\{v_n:n\in\omega\}\cup\{\textbf{c}_a:a\in A\}\cup\{\textbf{U},\textbf{R}\}.$$
Also let $$V(L_A)=\{v_n:n\in\omega\},$$
$$C(L_A)=\{\textbf{c}_a:a\in A\},$$
$$Term(L_A)=\{v_n:n\in\omega\}\cup\{\textbf{c}_a:a\in A\}$$
be  the classes of variables, constants and} terms of $L_A$, {\em respectively.}
\end{Def}

\begin{Def} \label{D:prima1}
{\em Given any set $A$, the class of} atomic formulas  {\em $AFml(L_A)$ of $L_A$ is defined as follows:
$$AFm(L_A)=\{\langle \equiv,t,s\rangle:t,s\in Term(L_A)\}\cup \{\langle \textbf{U},t\rangle:t\in Term(L_A)\}\cup$$$$\{\langle \textbf{R},t,s\rangle:t,s\in Term(L_A)\}.$$}
\end{Def}
[The meaning of the above codings is straightforward. $\langle \equiv,t,s\rangle$, $\langle\textbf{U},t\rangle$ and $\langle \textbf{R},t,s\rangle$ represent the formulas $t\equiv s$, $\textbf{U}(t)$ and $\textbf{R}(t,s)$, respectively.]

\begin{Def} \label{D:prima3}
{\em The class $Fml(L_A)$ of} formulas of $L_A$ {\em is defined as follows.
Define first the predicate $Fm_{L_A}(x,f,n)$ as follows:
$$Fml_{L_A}(x,f,n):=[n\neq 0 \wedge dom(f)=n \wedge  f(n-1)=x \ \wedge$$$$ (\forall k<n) [(f(k)\in AFml(L_A) \vee (\exists j<k)(f(k)=\langle \neg,f(j)\rangle) \ \vee$$
$$(\exists j,l<k)(f(k)=\langle \wedge,f(j),f(l)\rangle) \ \vee$$$$ (\exists j<k)(\exists m\in\omega)(f(k)=\langle \exists,m,f(j)\rangle)]].$$
Then we set
$$Fml(L_A)=\{x:(\exists f)(\exists n\in\omega)Fml_{L_A}(x,f,n)\}.$$
}
\end{Def}
The preceding definition of  formulas is essentially  the one  given in \cite[Def. 5.2 of Chapter 3]{Dr74}. It goes  smoothly despite the fact that  $\omega$ need not be a set. However the above  definition of $Fml(L_A)$ is not quite precise.  By writing $Fml_{L_A}(x,f,n)$, $f$ is intended to be  a function with domain $n$ that enumerates  the  subformulas of $x$ and {\em only them.}
But  if, for example,  $dom(f)=3$ and $f(0)$, $f(1)$ and $f(2)$ are atomic formulas, and $f(2)=x$, then $Fml(x,f,3)$ holds according to  \ref{D:prima3}, although  $f(0)$ and  $f(1)$ are not subformulas of $x$. Thus an additional constraint must be added to the definition of $Fml_{L_A}(x,f,n)$ in order to prevent $f$ from enumerating irrelevant atomic formulas. This is simply the requirement for the  domain of $f$  to be {\em minimal,} specifically that  $dom(f)=|Sub(x)|$ (the number  of  subformulas of $x$). In the above example the domain of a function enumerating the subformulas of an atomic formula should be $1$ not $3$.  This  requirement can be  formally expressed by a simple modification to  the definition of $Fml_{L_A}(x,f,n)$, and hence to that of $Fml(L_A)$.\footnote{The modification is this: We set
$$Fml(L_A)=\{x:(\exists f)(\exists n\in\omega)Fml^*_{L_A}(x,f,n)\},$$
where
$$Fml^*_{L_A}(x,f,n):=Fml_{L_A}(x,f,n) \wedge (\forall g)(\forall m)[Fml_{L_A}(x,g,m)\rightarrow n\leq m].$$
} Henceforth we  assume  that this requirement is implicitly satisfied whenever we write $Fml_{L_A}(x,f,n)$.

The crucial thing about  \ref{D:prima3}  is its capability to support inductive proofs and recursive definitions. We let the letters  $\sigma$, $\tau$, denote  elements of $Fml(L_A)$. As a first application of \ref{D:prima3}, every  $\sigma\in Fml(L_A)$ is assigned a {\em length}, which is the domain $n$ of some   enumerating function  $f$ for the subformulas of $\sigma$, or, since these functions are all minimal,   the number of subformulas $\sigma$.  More generally, in view of the validity of induction along  $\omega$ the following holds.

\begin{Lem} \label{L:formind}
Let $X\subseteq Fml(L_A)$ be a subclass of $Fml(L_A)$. If $AFml(L_A)\subseteq X$ and $X$ is closed with respect to $\neg$, $\wedge$ and $\exists$, then $X=Fml(L_A)$.
\end{Lem}

{\em Proof.} Assume $X$ is as stated and suppose $Fml(L_A)-X\neq \emptyset$.
Then, by the inductive properties of $\omega$, there is  $\sigma\in Fml(L_A)-X$  of least length $n$. Then we immediately obtain a contradiction from the definition of $Fml(L_A)$. \telos

\vskip 0.2in

Given a language $L_A$ as above, $L_A$-structures are defined as follows.

\begin{Def} \label{D:Lstructure}
{\em For any set $A$, an} $L_A$-structure {\em is a quadruple
$${\cal M}=\langle M,U,R,I\rangle,$$ where  $U\subseteq M$, $R\subseteq M\times M$, and $I$ is a (set) function $I:A\rightarrow M$. We refer to $M$ as the} domain {\em  of ${\cal M}$ and to $I$ as the} constant assignment {\em  for $L_A$. If either $\textbf{U}$ or $\textbf{R}$ is missing from $L_A$, the $L_A$-structures are triples ${\cal M}=\langle M,R,I\rangle$ or ${\cal M}=\langle M,U,I\rangle$, respectively. The} interpretation {\em of the extra-logical  symbols of  $L_A$ in  ${\cal M}$ is defined as follows: $\textbf{U}^{\cal M}=U$, $\textbf{R}^{\cal M}=R$, $\textbf{c}_a^{\cal M}=I(a)$, for each $a\in A$. On the other hand $\equiv^{\cal M}$ is the identity.}
\end{Def}

For every $\sigma\in Fml(L_A)$, the (finite) set of free variables of $\sigma$, denoted $FV(\sigma)$ is defined as usual by induction on the length of $\sigma$ (that is, along the steps of \ref{D:prima3}). Also for every $\sigma$, $|FV(\sigma)|\in\omega$. We come to the definition of the satisfaction relation $Sat({\cal M},\sigma, e)$ which formalizes the relation ${\cal M}\models \sigma(e(0),\ldots,e(n-1))$, for an $L_A$-structure ${\cal M}=\langle M,U,R,I\rangle$, a formula $\sigma(v_0,\ldots,v_{n-1})$ with free variables $v_0,\ldots,v_{n-1}$, and a mapping $e:n\rightarrow M$, that is, $e\in {}^nM$ (recall that by the discussion in the end of section 2.2, ${}^nM$ is a set). The next definition is an adaptation of Definition 5.4 of Chapter 3 of \cite{Dr74}. (Recall that the letters $f,g$ always denote functions.) In the definition below we assume that both symbols  \textbf{U} and \textbf{R} occur; if some of them is missing, the   definition is modified in the obvious way.

\begin{Def} \label{D:sat}
{\em Let  $Sat({\cal M},\sigma, e)$ denote the following relation:
$$(\exists M,U,R,I)(\exists f,g)(\exists n,m\in\omega)[{\cal M}=\langle M,U,R,I\rangle \ \wedge$$ $$Fml_{L_A}(\sigma,f,n)  \wedge |FV(\sigma)|=m \wedge dom(g)=n \wedge e\in g(n-1) \ \wedge $$
$$  (\forall k<n) S(k,f,g,m,{\cal M})],$$
where $S(k,f,g,m,{\cal M}):=$
$$(\exists i,j)[f(k)=\langle \equiv,v_i,v_j\rangle \wedge g(k)=\{h\in{}^mM:h(i)=h(j)\}] \ \vee$$
$$(\exists i)(\exists a\in A)[f(k)=\langle \equiv,v_i,\textbf{c}_a\rangle \wedge g(k)=\{h\in{}^mM:h(i)=I(a)\}] \ \vee$$
$$(\exists j)(\exists a\in A)[f(k)=\langle \equiv,\textbf{c}_a,v_j\rangle \wedge g(k)=\{h\in{}^mM:h(j)=I(a)\}] \ \vee$$
$$(\exists a\in A)[f(k)=\langle \equiv,\textbf{c}_a,\textbf{c}_a\rangle  \wedge g(k)={}^mM] \ \vee$$
$$(\exists a\neq b\in A)[f(k)=\langle \equiv,\textbf{c}_a,\textbf{c}_b\rangle  \wedge g(k)=\emptyset] \ \vee$$
$$(\exists i)[f(k)=\langle \textbf{U},v_i\rangle \wedge g(k)=\{h\in{}^mM:h(i)\in U\}] \ \vee$$
$$(\exists a\in A)[f(k)=\langle \textbf{U},\textbf{c}_a\rangle \wedge I(a)\in U \wedge  g(k)={}^mM] \ \vee$$
$$(\exists a\in A)[f(k)=\langle \textbf{U},\textbf{c}_a\rangle \wedge I(a)\notin U \wedge  g(k)=\emptyset] \ \vee$$
$$(\exists i,j)[f(k)=\langle \textbf{R},v_i,v_j\rangle \wedge g(k)=\{h\in{}^mM:\langle h(i),h(j)\rangle\in R\}] \ \vee$$
$$(\exists i)(\exists a\in A)[f(k)=\langle \textbf{R},v_i,\textbf{c}_a\rangle \wedge g(k)=\{h\in{}^mM:\langle h(i),I(a)\rangle\in R\}] \vee$$
$$(\exists j)(\exists a\in A)[f(k)=\langle \textbf{R},\textbf{c}_a,v_j\rangle \wedge g(k)=\{h\in{}^mM:\langle I(a),h(j)\rangle\in R\}] \ \vee$$
$$(\exists a,b\in A)[f(k)=\langle \textbf{R},\textbf{c}_a,\textbf{c}_b\rangle \wedge \langle I(a),I(b)\rangle\in R \wedge g(k)={}^mM] \ \vee$$
$$(\exists a,b\in A)[f(k)=\langle \textbf{R},\textbf{c}_a,\textbf{c}_b\rangle \wedge \langle I(a),I(b)\rangle\notin R \wedge g(k)=\emptyset] \ \vee$$
$$(\exists i)[f(k)=\langle \neg,f(i)\rangle \wedge g(k)={}^mM-g(i)] \ \vee$$
$$(\exists i,j)[f(k)=\langle \wedge,f(i),f(j)\rangle \wedge g(k)=g(i)\cap g(j)] \ \vee$$
$$(\exists i,j)[f(k)=\langle \exists,i,f(j)\rangle \wedge \ g(k)=\{h\in {}^mM:(\exists x\in M) (h(i/x)\in g(j))\}],$$
(where in the last clause $h(i/x)$ is the tuple of ${}^mM$ resulting from $h$ if we replace $\langle i,h(i)\rangle$ with $\langle i,x\rangle$).}
\end{Def}

In the preceding definition the function $g$, that enumerates the sets of  assignments that make true  the subformulas of $\sigma$, is definable in EST because  on the one hand ${}^mM$ is a set, and on the other hand in each  clause of the definition,  $g(k)$ is a $\Delta_0$ subclass of ${}^mM$, therefore a set.

Having defined what ${\cal M}\models\sigma(x_0,\ldots,x_{m-1})$ means for an $L_A$-structure ${\cal M}$, an $L_A$-formula $\sigma(v_0,\ldots,$ $v_{m-1})$  with its free variables being among $v_0,\ldots,v_{m-1}$, and for $x_0,\ldots, x_{m-1}\in M$, we can then define  elementary embeddings  from one  $L_A$-structure into another as usual.

Given two $L_A$-structures ${\cal M}$,  ${\cal N}$, we say that ${\cal M}$ is {\em elementarily embeddable} in ${\cal N}$, notation ${\cal M}\precsim {\cal N}$, if there is a 1-1 function $f:M\rightarrow N$ such that for every formula $\sigma(v_0,\ldots,v_{m-1})$ of $L_A$ with free variables among $v_0,\ldots,v_{m-1}$, and any  $x_0,\ldots, x_{m-1}\in M$,
$${\cal M}\models\sigma(x_0,\ldots,x_{m-1})\leftrightarrow {\cal N}\models\sigma(f(x_0),\ldots,f(x_{m-1})).$$
Let $f:{\cal M}\precsim {\cal N}$ denote the fact that $f$ is an elementary embedding of ${\cal M}$ into  ${\cal N}$. Sometimes, for more precision,  we need to specify the language we refer to. Then we say that $f:M\rightarrow N$ is an {\em $L_A$-elementary embedding} and we  denote it by
$$f:{\cal M}\precsim_{L_A} {\cal N}.$$
Clearly the last relation  is definable in EST.
The following simple fact will be repeatedly used below.
\begin{Fac} \label{F:elem}
Let ${\cal M}=\langle M,U,R,I\rangle$, ${\cal N}=\langle N,Z,S,J\rangle$ be  $L_A$-structures and $f: M\rightarrow N$ be an $L_A$-elementary embedding (or  just an $L_A$-embedding).  Then $f\circ I=J$, that is, for every $a\in A$, $f(I(a))=J(a)$. Equivalently, for every $a\in A$, $f(\textbf{c}_a^{\cal M})=\textbf{c}_a^{\cal N}$.
\end{Fac}

{\em Proof.} For every $L_A$-embedding $f:{\cal M}\rightarrow {\cal N}$, by definition $f(\textbf{c}_a^{\cal M})=\textbf{c}_a^{\cal N}$. Also for every  $a\in A$,  $I(a)=\textbf{c}_a^{\cal M}$ and $J(a)=\textbf{c}_a^{\cal N}$. Therefore $f(I(a))=J(a)$.  \telos

\vskip 0.2in

The language ${\cal L}=\{\in\}$ of EST is just a particular instance of the languages $L_A$ above, namely one with $A=\emptyset$ and one binary relation symbol $\in$.  So the classes of formulas and  sentences of ${\cal L}$ are already definable in EST.

Given a formula $\phi(x)$ of ${\cal L}=\{\in\}$, let  $X_\phi$ denote the extension  of $\phi$, although $X_\phi$ needs not be a set  in EST. The expression ``$X_\phi$ is a proper class''  simply stands for  the   ${\cal L}$-sentence:
$$(\forall x)(\exists y)[(\phi(y)\wedge y\notin x) \vee (\neg\phi(y)\wedge y\in x)].$$
\noindent [Note that in ZF ``$X_\phi$ is a proper class'' is formulated  just as $(\forall x)(X_\phi\not\subseteq x)$ because of Separation. But in EST, where Separation  is missing and $X_\phi\subseteq x$ does not imply that $X_\phi$ is a set, ``$X_\phi$ is a proper class'' need to be formulated  as $(\forall x)(X_\phi\neq x)$.]

Given any formula $\phi(x)$ of ${\cal L}$ in one free variable, the instance ${\rm VP}_\phi$ of Vop\v{e}nka's Principle is formulated as follows:

\vskip 0.1in

${\rm VP}_\phi$: {\em ``For every set $A$,  if $X_\phi$ is a proper class of  $L_A$-structures,   then there are ${\cal M}\neq {\cal N}\in X_\phi$ and $f$ such that  $f:{\cal M}\precsim_{L_A}{\cal N}$''.}

\vskip 0.1in

Clearly ${\rm VP}_\phi$ is an ${\cal L}$-sentence so the class
$${\rm VP}=\{x: (\exists\phi\in Fml({\cal L}))(x={\rm VP}_\phi)\}$$
is definable in EST. This completes the description of the scheme VP in EST.

\subsection{Consequences of EST+VP}

We begin  with the proof of Infinity because  it does not depend on any form of Replacement or Powerset. We show that the class $\omega=\{x:Nat(x)\}$, as defined in section 2.2, is a set in EST+VP.

\begin{Thm} \label{T:inf}
In ${\rm EST+VP}$, the class $\omega=\{x:Nat(x)\}$ is a set. Therefore ${\rm EST+VP}\vdash Inf$.
\end{Thm}

{\em Proof.} To avoid dealing with elementary embeddings of $\emptyset$, let $\omega^*=\omega-\{0\}$. Obviously in EST $\omega$ is a set iff  $\omega^*$ is a set. So towards reaching a contradiction assume that $\omega^*$ is a proper class. Consider the language $L_\emptyset=\{\textbf{R}\}$ with only a binary relation symbol $\textbf{R}$ (that is, $L_\emptyset$ contains no  constants $\textbf{c}_a$). For each $x\in\omega^*$, let $\in_x=\{\langle y,z\rangle\in x\times x:y\in z\}$ (the restriction of $\in$ to $x$). Since $x\times x$ is a set, by $\Delta_0$-Separation $\in_x$  is a set too. Let ${\cal M}_x=\langle x,\in_x\rangle$. Each ${\cal M}_x$ is an $L_\emptyset$-structure by interpreting $\textbf{R}$ by $\in_x$. Let also $$K=\{{\cal M}_x:x\in \omega^*\}.$$
We claim that $K$ is a proper class when $\omega^*$ is so. Indeed, assume $K$ is a set. Since
$$K=\{\langle x,\in_x\rangle:x\in \omega^*\}=\{\{\{x\},\{x,\in_x\}\}:x\in \omega^*\},$$
clearly  $\omega^*\subset \cup^2K=\cup\cup K$. In particular $\omega^*=\{x\in \cup^2K:Nat(x)\}$. Since $\cup^2K$ is a set and $Nat(x)$ is $\Delta_0$, $\omega^*$ is a set by $\Delta_0$-Separation, a contradiction.

Thus $K$ is a proper class of $L_\emptyset$-structures. By VP there exist  $x_0\neq x_1\in \omega^*$ and a function $f:x_0\rightarrow x_1$ such that  $f:{\cal M}_{x_0}\precsim {\cal M}_{x_1}$.

But one can easily see by the  clauses of Fact \ref{F:Natfacts} that this is false.
Indeed, since $f$ is 1-1 and  $x_0\neq x_1$, by \ref{F:Natfacts} (i),(ii), (iii),   $x_0\in x_1$, and hence $x_0$ is a proper initial segment of $x_1$.  Also, by elementarity, we can see by induction on $x_0$ that $f$ is the identity. [$0=\emptyset$ is the first element of both $x_0$ and $x_1$, so $f(0)=0$. Inductively, if $f(n)=n$ then the next element of $n$ should be sent to the next element of $f(n)$, that is, $f(n+1)=n+1$.] The elements of $\omega^*$ are all successor ordinals, so let  $x_0=y_0\cup\{y_0\}$, $x_1=y_1\cup\{y_1\}$. Since $x_0\neq x_1$, we have $y_0\neq y_1$. By \ref{F:Natfacts} (iv),  $y_0$ is the greatest element of $x_0$ and $f(y_0)=y_0$, since $f$ is the identity, while, by elementarity,  $f(y_0)$ should be the greatest  element $y_1$ of $x_1$. But $y_1\neq y_0$, a contradiction. \telos

\vskip 0.2in

Next we come to the proof of Replacement.

\begin{Thm} \label{T:Deltanew}
(i)  ${\rm EST+VP}\vdash \Delta_0$-Rep.

(ii) ${\rm EST+VP}+\Delta_0$-$Rep\vdash$ Rep.

(iii)  Therefore ${\rm EST+VP}\vdash$ Rep.
\end{Thm}

{\em Proof.} (i) To prove $\Delta_0$-$Rep$, let $\phi(x,y)$ be a $\Delta_0$ formula such that $(\forall x)(\exists !y)$ $\phi(x,y)$. This  defines a class mapping $F_\phi:V\rightarrow V$ such that $F_\phi(x)=y$ iff $\phi(x,y)$. Fix a set $A$. It suffices to show that the class $B=F_\phi``A=\{F_\phi(a):a\in A\}$ is a set. Let $L_A=\{\textbf{U}\}\cup\{\textbf{c}_a:a\in A\}$ be the language with  a unary relation symbol $\textbf{U}$ and a constant $\textbf{c}_a$ for each $a\in A$. For every $b\in B$, $A\times\{b\}$ is a set, so for each such $b$ consider the  $L_A$-structure $${\cal M}_b=\langle A\times\{b\},U_b, I_b\rangle,$$ where
 $U_b\subseteq A\times \{b\}$ is defined as follows: For every $a\in A$
$$\langle a,b\rangle\in U_b\iff F_\phi(a)=b.$$
We have  $U_b=\{\langle a,b\rangle\in A\times \{b\}:\phi(a,b)\}$. By $CartProd$, $A\times \{b\}$ is a set and since $\phi$ is $\Delta_0$,  $U_b$ is a set, by $\Delta_0$-$Sep$, that  interprets $\textbf{U}$, that is,  $\textbf{U}^{{\cal M}_b}=U_b$. The constant assignment $I_b:A\rightarrow A\times \{b\}$ is defined by  $I_b(a)=\langle a,b\rangle$, for each $a\in A$.
$I_b$ is a set too, by $\Delta_0$-$Sep$, because  $I_b=\{\langle x,\langle y,b\rangle\rangle\in A\times (A\times \{b\}):x=y\}$ and $x=y$ is $\Delta_0$.
This means that  for all $a\in A$ and $b\in B$,   $\textbf{c}_a^{{\cal M}_b}=\langle a,b\rangle$.   Let $$S=\{{\cal M}_b:b\in B\}.$$ It suffices to show that $S$ is a set. For suppose  that this is the case. Then clearly  for some $n\in\omega$ (actually for $n=7$), $B\subset \cup^nS$. Moreover
$$B=\{y\in \cup^nS:(\exists x\in A)(F_\phi(x)=y)\}=\{y\in \cup^nS:(\exists x\in A)\phi(x,y)\}.$$
Since $\cup^nS$ is a set and the formula $(\exists x\in A)\phi(x,y)$ is $\Delta_0$, it follows by  $\Delta_0$-Separation that $B$ is  set.

So let us verify that $S$ is a set. To reach a contradiction assume that $S$ is a proper class.  Then by VP there are $b,c\in B$, $b\neq c$, and a mapping $f:A\times\{b\}\rightarrow A\times\{b\}$ such that  $f:{\cal M}_b\precsim {\cal M}_c$. By elementarity, for every $a\in A$,
$$f(\langle a,b\rangle)=f(\textbf{c}_a^{{\cal M}_b})=\textbf{c}_a^{{\cal M}_c}=\langle a,c\rangle.$$ On the other hand, by elementarity again, for every $a\in A$, $$F_\phi(a)=b\Leftrightarrow \langle a,b\rangle\in U_b\Leftrightarrow f(\langle a,b\rangle)\in U_c\Leftrightarrow \langle a,c\rangle\in U_c\Leftrightarrow F_\phi(a)=c,$$
which is a contradiction since $b\neq c$.

(ii)  Now we work in EST+VP+$\Delta_0$-$Rep$, and prove that full Replacement holds. The proof is for the most part similar to that of clause (i) above. Let $\phi(x,y)$ be a formula such that $(\forall x)(\exists !y)\phi(x,y)$, and let $F_\phi(x)=y$ iff $\phi(x,y)$. We fix again a set $A$ and show that if $B=F_\phi``A$, then $B$ is a set. We define the structures ${\cal M}_b$ as before and we set  $S=\{{\cal M}_b:b\in B\}$. As in (i), it follows by means of VP that  $S$ cannot be a proper class. Thus $S$ is a set. The only departure from the proof of (i) is at the point of inferring that $B$ is a set from  $S$ being a set. This now can be inferred  by he help of $\Delta_0$-Replacement: just observe that the  mapping $S\ni {\cal M}_b\mapsto b\in B$ is clearly  $\Delta_0$-definable and onto. Therefore $B$ is a set.

(iii) Immediate from (i) and (ii). \telos

\vskip 0.2in

Now we come to the proof of Powerset, which is based on a clause of Theorem \ref{T:Deltanew}.

\begin{Thm} \label{T:Pow}
${\rm EST+VP}+\Delta_0$-Rep $\vdash Pow$. Therefore, by Theorem  \ref{T:Deltanew} (i), ${\rm EST+VP}\vdash Pow$.
\end{Thm}

{\em Proof.} We work in ${\rm EST+VP}+\Delta_0$-$Rep$. Fix a set $A$ and let again $L_A=\{\textbf{U}\}\cup\{\textbf{c}_a:a\in A\}$, where $\textbf{U}$ is  a unary relation symbol. For each $X\in {\cal P}(A)$ consider the $L_A$-structure $${\cal M}_X=\langle A,X,id_A\rangle,$$  where for each $X\subseteq A$, $\textbf{U}^{{\cal M}_X}=X$ and the constant assignment is the identity mapping $id_A:A\rightarrow A$. Note that $id_A$ is a set in EST,  by $\Delta_0$-$Sep$, since $id_A=\{\langle x,y\rangle\in A\times A:x=y\}$. Thus $\textbf{c}_a^{{\cal M}_X}=a$ for every $a\in A$.
To reach a contradiction, assume that ${\cal P}(A)$ is a proper class. Let
$$K=\{{\cal M}_X:X\in {\cal P}(A)\}.$$
The mapping $K\ni {\cal M}_X\mapsto X\in {\cal P}(A)$ is clearly $\Delta_0$, so by $\Delta_0$-Replacement, the class $K$  is proper too. By VP there are $X\neq Y\in {\cal P}(A)$ and an elementary embedding $f:{\cal M}_X\rightarrow {\cal M}_Y$. But then $f(a)=f(\textbf{c}_a^{{\cal M}_X})=\textbf{c}_a^{{\cal M}_Y}=a$, for every $a\in A$. That is, $f=id_A$. On the other hand, by elementarity, $f$ should map 1-1 $X$ onto $Y$, hence $X=Y$,  a contradiction. \telos

\vskip 0.2in

Now ${\rm ZF}={\rm EST}$+$\{Inf, Pow, Rep, Found\}$. So from Theorems \ref{T:inf}, \ref{T:Deltanew} and \ref{T:Pow}  we obtain immediately the following:

\begin{Cor} \label{C:mainnew}
${\rm EST}+Found+{\rm VP}={\rm ZF+VP}$, and ${\rm EST}+Found+{\rm AC+VP}={\rm ZFC+VP}$.
\end{Cor}

As  mentioned in the beginning of this section, the replacement of ${\rm VP}_4$ (used in Theorem \ref{T:main}) by VP (used in Theorems  \ref{T:inf},  \ref{T:Deltanew} and \ref{T:Pow}) was necessitated by the fact that  theorem V-P-H was  among the  assumptions of \ref{T:main}, while this is not the case for  \ref{T:inf}, \ref{T:Deltanew} and \ref{T:Pow}. Inspecting the proof of V-P-H in \cite{VPH65}, we see that it relies heavily on AC, as well as on the following two facts: (a) For every well ordered set $\langle x,\leq\rangle$, there is a (unique) ordinal $\alpha$ such that $\langle x,\leq\rangle\cong \langle \alpha,\in\rangle$. (b) For every ordinal $\alpha$, there exists the set of ordinals of countable cofinality  below $\alpha$, $\{\beta<\alpha:{\rm cf}(\beta)=\omega\}$. (a) requires $\Delta_1$-Replacement, while  (b)   requires, firstly, that $\omega$ is a set and, secondly, $\Sigma_1$-Separation.  These being available the proof of V-P-H goes through, so
$${\rm EST}+\{{\rm AC},``\omega \ \mbox{is a set''},\Delta_1\mbox{-}Rep, \Sigma_1\mbox{-}Sep\}\vdash \mbox{V-P-H}.$$
With  V-P-H at hand we can work exactly as in the proof of \ref{T:main} with ${\rm VP}_4$ in place of  VP  (a rigid binary relation $R$ on any set $A$ does the job that the constants  $\textbf{c}_a$ do in VP). Thus  we  obtain the following.

\begin{Cor} \label{C:Pultr}
The theory $${\rm EST}+\{{\rm AC},``\omega \ \mbox{is a set''},\Delta_1\mbox{-}Rep, \Sigma_1\mbox{-}Sep\}+{\rm VP}_4$$ proves Replacement and Powerset.
\end{Cor}

Corollary \ref{C:Pultr} is in sharp contrast to the results \ref{T:Deltanew} and \ref{T:Pow} above, which together show that EST+VP {\em alone} proves Replacement and Powerset. This gives a measure of the difference in apparent strength between principles ${\rm VP}_4$ and VP over EST.

\begin{Rem} \label{R:reduce}
{\em It is further worth mentioning that in the proofs of  theorems \ref{T:Deltanew} and \ref{T:Pow} we did not use  the full strength of VP. A simple inspection of the proofs shows that in these results we used the fact that for a given proper class $X_\phi$ of $L_A$-structures there are distinct structures ${\cal M}$, ${\cal N}$ in $X_\phi$ and an} embedding only $f:{\cal M}\rightarrow {\cal N}$, {\em  rather than an elementary embedding. In contrast, in Theorem \ref{T:inf} some kind of elementarity for $f$ is required. Therefore   \ref{T:Deltanew} and \ref{T:Pow} can still be established by means of the following weaker form, ${\rm VP}_0$, of Vop\v{e}nka's Principle. Let ${\rm VP}_0$ result from VP if  ``elementary embedding'' is replaced by ``embedding'' (while the languages $L_A$ still contain finitely many relations and an arbitrary set of constants). Then ${\rm EST}+{\rm VP}_0$ proves Replacement and Powerset.}
\end{Rem}

Since ${\rm EST}\subset {\rm LZFC}$ (see footnote 1, after Theorem \ref{T:main}), as an immediate Corollary to Theorems \ref{T:Deltanew} and \ref{T:Pow} we obtain the  following improvement to  Theorem 6.3 of \cite{Tz14}:

\begin{Thm} \label{T:mainnew1}
${\rm LZFC+VP}$ proves Replacement and Powerset, that is, ${\rm ZFC}$ $\subseteq {\rm LZFC+VP}$.
\end{Thm}

The  improvement consists of course in  the fact that the requirement for $T$  to prove  V-P-H is no longer needed for LZFC.

\vskip 0.2in

On the other hand, ZFC+VP implies the existence of a proper class of extendible  cardinals (see Lemma 20.25 of \cite{Je03} and the remark immediately after its  proof). Since for every extendible cardinal $\kappa$, $V_\kappa\models{\rm ZFC}$, we immediately infer that ${\rm ZFC+VP}\vdash Loc({\rm ZFC})$, where $Loc({\rm ZFC})$ is the central axiom of LZFC saying that ``every set belongs to a transitive model of ZFC.'' From this we have:

\begin{Lem} \label{L:jech}
  ${\rm LZFC}\subset {\rm ZFC+VP}$.
\end{Lem}
From   Theorem \ref{T:mainnew1} and Lemma \ref{L:jech} we  obtain:

\begin{Thm} \label{T:main2}
 ${\rm LZFC+VP}={\rm ZFC+VP}$.
 \end{Thm}

\section{VP and Foundation}
Let ${\rm ZF}_0={\rm ZF}-\{Found\}$. In this subsection we show that VP does not prove Foundation over  ${\rm ZF}_0$. Namely, the following holds:

\begin{Thm} \label{T:Found}
If ${\rm ZF+VP}$ is consistent, then so is ${\rm ZF}_0+{\rm VP}+\neg Found$. Similarly with {\rm ZFC} in place of {\rm ZF}.
\end{Thm}

The proof is by the well-known method of using  a nonstandard membership relation $\in_\pi$ in $V$, produced  by a definable permutation of $V$. Namely, it is a rather folklore result that if $V$ is the universe of ZFC,  $\pi:V\rightarrow V$ is a definable permutation, and $\in_\pi$ is the binary relation defined by $x\in_\pi y$ iff $x\in \pi(y)$, then $V_\pi\models {\rm ZFC}_0$ (see \cite[Ch. IV, exercise 18]{Ku83}). For simplicity, let us abbreviate henceforth $\langle V,\in_\pi\rangle$ by $V_\pi$ and $\langle V,\in\rangle$ by $V$. In order for Foundation to fail in $V_\pi$ it suffices to take $\pi$ so that $x\in \pi(x)$ for some $x$.  In this way we shall prove the next Theorem from which Theorem \ref{T:Found} follows.

\begin{Thm} \label{T:VPreserve}
If $V\models {\rm ZF+VP}$, then there is a permutation $\pi:V\rightarrow V$ such that  $V_\pi\models {\rm VP}+\neg Found$.
\end{Thm}

We shall need first some preliminary  definitions and Lemmas. Given a  definable permutation $\pi:V\rightarrow V$, let us denote by $\phi^\pi$, for every formula $\phi$ of ${\cal L}=\{\in\}$, the formula resulting from $\phi$ if we replace every atomic subformula $x\in y$ occurring in $\phi$ by $x\in \pi(y)$. The following are easy to check by induction on the length of $\phi$.

\begin{Lem} \label{L:Prel}
(i) The mapping $\phi\mapsto \phi^\pi$ commutes with connectives and quantifiers, that is, $(\phi\rightarrow \psi)^\pi=(\phi^\pi\rightarrow \psi^\pi)$, $(\neg\phi)^\pi=\neg\phi^\pi$, and $(\forall x\phi)^\pi=(\forall x)\phi^\pi$.

(ii) For every sentence $\phi$, $V_\pi\models\phi$ iff $V\models\phi^\pi$.
\end{Lem}

For  every standard notion of $V$, like  singleton, ordered pair, $n$-tuple, relation, function, there is a corresponding $\pi$-notion for $V_\pi$. For instance a $\pi$-pair is a set $z$ such that $V_\pi\models$ ``$z$ is an ordered pair''. The latter holds iff for some  $x,y$ $V_\pi\models z=\langle x,y\rangle$. By \ref{L:Prel} (ii), this is equivalent to $(z=\langle x,y\rangle)^\pi$ (we often  write just  $\phi$ instead of $V\models\phi$). Also the latter is more conveniently denoted by $z=\langle x,y\rangle^\pi$. Similarly the fact that $Q$ is a $\pi$-binary relation between sets  $M$ and $N$ means that $(Q\subseteq M\times N)^\pi$ is true. A $\pi$-notion expressed by a sentence $\phi$ is said to be {\em absolute} if $\phi^\pi\leftrightarrow \phi$. In the next two Lemmas we give some simple sufficient conditions concerning the permutation $\pi$ in order for some key notions to be absolute.

\begin{Lem} \label{L:key1}
Suppose that  $\pi:V\rightarrow V$ fixes all finite sets. Then:

(i)  $\pi$-pairs are absolute, that is, for all $x,y$, $\langle x,y\rangle^\pi=\langle x,y\rangle$.

(ii) For all $M,N,Q$,  $(Q\subseteq M\times N)^\pi\leftrightarrow \pi(Q)\subseteq \pi(M)\times \pi(N)$.

(iii) $[f:M\rightarrow N \ \mbox{is a function}]^\pi$ is equivalent to $\pi(f):\pi(M)\rightarrow \pi(N) \ \mbox{is a function}$.

\end{Lem}

{\em Proof.} (i)  Assume $\pi$ fixes all finite sets. Then so does also $\pi^{-1}$.  Analyzing the definition of  $(z=\langle x,y\rangle)^\pi:=(z=\{\{x\},\{x,y\}\})^\pi$, it is easy to see that
\begin{equation} \label{E:pi}
z=\langle x,y\rangle^\pi \leftrightarrow z=\pi^{-1}(\{\pi^{-1}(\{x\}),\pi^{-1}(\{x,y\})\}).
\end{equation}
All the arguments of $\pi^{-1}$ in the right-hand side of (\ref{E:pi}) are finite, so $\pi^{-1}$ fixes them. Therefore (\ref{E:pi}) implies
$$z=\langle x,y\rangle^\pi \leftrightarrow z=\{\{x\},\{x,y\}\}=\langle x,y\rangle.$$

(ii) Analyzing the definition of $(Q\subseteq M\times N)^\pi$ we  see that
\begin{equation} \label{E:pi-binary}
(Q\subseteq M\times N)^\pi \leftrightarrow \pi(Q)\subseteq \{\langle x,y\rangle^\pi: x\in \pi(M),y\in \pi(N)\}.
\end{equation}
By  (i) above, $\langle x,y\rangle^\pi=\langle x,y\rangle$ for every pair. So (\ref{E:pi-binary}) gives
$$(Q\subseteq M\times N)^\pi \leftrightarrow \pi(Q)\subseteq \{\langle x,y\rangle: x\in \pi(M),y\in \pi(N)\}=\pi(M)\times\pi(N).$$

(iii) Easy to check. \telos

\vskip 0.2in

Recall that a language $L_A$ consists of the symbols  $\textbf{R}$,  $\textbf{U}$ and  $\textbf{c}_a$, for $a\in A$, and an $L_A$-structure is a quadruple ${\cal M}=\langle M,U,R,I\rangle$, where $R\subseteq M\times M$, $U\subseteq M$ and $I:A\rightarrow M$ is a mapping.

Given an $L_A$ structure ${\cal M}=\langle M,U,R,I\rangle$ and a permutation $\pi:V\rightarrow V$, let us set
$${\cal M}^\pi:=\langle \pi(M),\pi(U),\pi(R),\pi(I)\rangle.$$

\begin{Lem} \label{L:I-function}
Suppose that $\pi:V\rightarrow V$ fixes all finite sets and $\omega$.

(i) Let $L_A$ be a first-order language  in the sense of $V_\pi$, for some $A\in V_\pi$.  Then $\pi(L_A)=L_{\pi(A)}$ is a language in $V$.

(ii) If  $\sigma$ is a formula of $L_A$ in  the sense of  $V_\pi$, then $\sigma$ is a formula of $L_{\pi(A)}$ in $V$.

(iii) If ${\cal M}$ is an  $L_A$-structure in $V_\pi$, then  ${\cal M}^\pi$ is an  $L_{\pi(A)}$-structure in $V$.

(iv) If ${\cal M}$ is an  $L_A$-structure in $V_\pi$, then
$$V_\pi\models [x_0,x_1\in M \wedge ({\cal M}\models {\rm \textbf{R}}(x_0,x_1))]\leftrightarrow V\models [x_0,x_1\in \pi(M) \wedge ({\cal M}^\pi\models {\rm \textbf{R}}(x_0,x_1))],$$
and similarly for the predicate $\textbf{U}$.
\end{Lem}

{\em Proof.} (i) Recall that by Definition \ref{D:prima},
$$V_\pi\models L_A=\{2,3,4,5,6,7\}\cup\{\langle 0,n\rangle:n \in \omega\}\cup\{\langle 1,a\rangle:a\in A\}.$$ In view of Lemma \ref{L:Prel}, this is equivalently written
$$V\models\pi(L_A)=\{2,3,4,5,6,7\}^\pi\cup\{\langle 0,n\rangle^\pi:n \in \pi(\omega)\}\cup\{\langle 1,a\rangle^\pi:a\in \pi(A)\}.$$
Since $\pi$ fixes all finite sets and $\omega$, $\pi$-pairs are absolute by Lemma \ref{L:key1},  $\{2,3,4,5,6,7\}^\pi=\{2,3,4,5,6,7\}$ and $\pi(\omega)=\omega$, so
$$V\models\pi(L_A)=\{2,3,4,5,6,7\}\cup\{\langle 0,n\rangle:n \in \pi(\omega)\}\cup\{\langle 1,a\rangle:a\in \pi(A)\}.$$
But  the right-hand side of the above equation is clearly the language $L_{\pi(A)}$, so $\pi(L_A)=L_{\pi(A)}$.

(ii) Since $\pi(\omega)=\omega$ and $\pi(n)=n$ for every $n\in\omega$, $\omega$ is absolute in $V_\pi$. Thus the claim follows by a simple induction on   $\sigma$ along the steps of Definition \ref{D:prima3}.

(iii) That ${\cal M}=\langle M,U,R,I\rangle$ is an $L_A$-structure in the sense of $V_\pi$ means that $(U\subseteq M)^\pi$, $(R\subseteq M\times M)^\pi$,  and $(I:A\rightarrow M$ is a function)$^\pi$.  Since $\pi$ fixes all finite sets, these facts are translated into  $V$, according to \ref{L:key1}, as $\pi(U)\subseteq \pi(M)$,   $\pi(R)\subseteq \pi(M)\times \pi(M)$ and $\pi(I):\pi(A)\rightarrow \pi(M)$ is a function, respectively. But this means that $\langle \pi(M),\pi(U),\pi(R),\pi(I)\rangle$, that is, ${\cal M}^\pi$, is an $L_{\pi(A)}$-structure.

(iv) Let ${\cal M}=\langle M,U,R,I\rangle$ be an $L_A$-structure in $V_\pi$ and let $$V_\pi\models [x_0,x_1\in M \wedge ({\cal M}\models {\rm \textbf{R}}(x_0,x_1))].$$ Obviously this is  equivalently written $$V_\pi\models [x_0,x_1\in M \wedge \langle x_0,x_1\rangle \in R)].$$ Its translation to $V$ is $$V\models [x_0,x_1\in \pi(M) \wedge \langle x_0,x_1\rangle^\pi \in \pi(R))].$$ Since $\langle x_0,x_1\rangle^\pi=\langle x_0,x_1\rangle$ by the condition on $\pi$, the latter also is equivalent to $$V\models [x_0,x_1\in \pi(M) \wedge ({\cal M}^\pi\models {\rm \textbf{R}}(x_0,x_1))].$$ \telos

\begin{Lem} \label{L:key2}
Let  $\pi$ be a permutation that fixes all finite sets and $\omega$, and let  ${\cal M}=\langle M,U,R,I\rangle$, ${\cal N}=\langle N,Z,S,J\rangle$ be $L_A$-structures in $V_\pi$. Then for any   $f$,
$$V_\pi\models [f:{\cal M}\precsim_{L_A}{\cal N}]$$ iff
$$V\models [\pi(f):{\cal M}^\pi\precsim_{L_{\pi(A)}}{\cal N}^\pi].$$
\end{Lem}

{\em Proof.} Let us sketch the proof of direction ``$\rightarrow$''. The other direction is similar. Assume that $\pi$ is as stated,  ${\cal M}$, ${\cal N}$ are $L_A$-structures in $V_\pi$ and $V_\pi\models [f:{\cal M}\precsim_{L_A}{\cal N}]$, that is, $f:M\rightarrow N$ is an $L_A$-elementary embedding.  By Lemma \ref{L:I-function}, ${\cal M}^\pi$, ${\cal N}^\pi$ are $L_{\pi(A)}$-structures in $V$. We have to show that $\pi(f):\pi(M)\rightarrow \pi(N)$ is an $L_{\pi(A)}$-elementary embedding.
We must show that  for every $L_{\pi(A)}$ formula $\sigma(v_0,\ldots,v_{n-1})$ and any $x_0,\ldots,x_{n-1}\in \pi(M)$,
$${\cal M}^\pi\models\sigma(x_0,\ldots,x_{n-1})\leftrightarrow {\cal N}^\pi\models\sigma(\pi(f)(x_0),\ldots,\pi(f)(x_{n-1})).$$ This is shown by routine induction on the length of $\sigma$. Let us just show  the above for the  atomic sentences ${\bf R}(x_0,x_1)$ of $L_{\pi(A)}$.   This amounts to showing that if  $x_0,x_1\in \pi(M)$ and $y_0,y_1\in \pi(N)$, then (in $V$):
$$\langle x_0,y_0\rangle\in \pi(f) \wedge \langle x_1,y_1\rangle\in \pi(f) \rightarrow$$
\begin{equation} \label{E:phi5}
({\cal M}^\pi\models {\bf R}(x_0,x_1)\leftrightarrow {\cal N}^\pi\models {\bf R}(y_0,y_1)).
\end{equation}
But by our assumption $V_\pi\models[f:{\cal M}\precsim_{L_A}{\cal N}]$, we have that for all $x_0,x_1\in M$ and $y_0,y_1\in N$:
$$V_\pi\models[\langle x_0,y_0\rangle\in f \wedge \langle x_1,y_1\rangle\in f \rightarrow$$
\begin{equation} \label{E:phip5}
({\cal M}\models {\bf R}(x_0,x_1)\leftrightarrow {\cal N}\models {\bf R}(y_0,y_1)).
\end{equation}
By Lemmas \ref{L:Prel} and \ref{L:I-function} (iv) (that holds because of our conditions  about $\pi$), (\ref{E:phi5}) is just the translation of (\ref{E:phip5}) to $V$.
The other steps of the induction are routine. This completes the proof of the Lemma.\telos

\vskip 0.2in

Now we are ready to prove \ref{T:VPreserve}.

\vskip 0.2in

{\em Proof of Theorem \ref{T:VPreserve}.}  Let $V\models {\rm ZF+VP}$. Pick a permutation $\pi$ of $V$ that  fixes all finite sets and $\omega$. Then Lemmas  \ref{L:key1}, \ref{L:I-function} and  \ref{L:key2} above hold. Suppose also that for some  sets $X$, $Y$ such that $X\in Y$, $\pi$ exchanges $X$ and $Y$, that is,  $\pi(X)=Y$, $\pi(Y)=X$,  so Foundation fails in $V_\pi$. It remains to show that $V_\pi\models {\rm VP}$.

Let $L_A$ be a language in the sense of $V_\pi$. The set $A$, that essentially contains the constants of the language, can be arbitrary, so in particular we may have $A=X$ or $A=Y$. For that reason in general $\pi(A)\neq A$.   Let $\phi(x)$ be a formula of ${\cal L}=\{\in\}$ such that
\begin{equation} \label{E:assum}
V_\pi\models \mbox{ ``$X_\phi$ is a proper class of $L_A$-structures''}.
\end{equation}
We have to  show that
\begin{equation} \label{E:concl}
V_\pi\models (\exists{\cal M}\neq {\cal N}\in X_\phi)({\cal M}\precsim_{L_A} {\cal N}).
\end{equation}
Now (\ref{E:assum}) implies that  $X_\phi$ is a proper class in $V_\pi$, that is,
\begin{equation} \label{E:proper1}
V_\pi\models (\forall x)(\exists y)(\phi(y) \wedge y\notin x),
\end{equation}
and also
\begin{equation} \label{E:omited1}
V_\pi\models (\forall x)[\phi(x)\rightarrow \ \mbox{$x$ is an $L_A$-structure}).
\end{equation}
From (\ref{E:proper1}) we have that
\begin{equation} \label{E:proper2}
V\models (\forall x)(\exists y)(\phi^\pi(y) \wedge y\notin \pi(x)),
\end{equation}
that is, $X_{\phi^\pi}$ is a proper class in $V$. Moreover if ${\cal M}\in X_{\phi^\pi}$, then $V_\pi\models \phi({\cal M})$, so by (\ref{E:omited1})  ${\cal M}$ is an $L_A$-structure in $V_\pi$. By Lemma \ref{L:I-function} (iii), ${\cal M}^\pi$ is an $L_{\pi(A)}$-structure in $V$. Thus
\begin{equation} \label{E:xeka}
{\cal M}\in X_{\phi^\pi}\rightarrow {\cal M}^\pi \ \mbox{is an $L_{\pi(A)}$-structure}.
\end{equation}
Consider the formula $\psi(x)$ of ${\cal L}$ defined by
$$\psi(x):=(\exists {\cal M})(x={\cal M}^\pi \wedge \phi^\pi({\cal M})).$$
Then clearly for every ${\cal M}$
\begin{equation} \label{E:star}
V\models\phi^\pi({\cal M})\iff V_\pi\models \phi({\cal M})\iff V\models \psi({\cal M}^\pi).
\end{equation}
By (\ref{E:xeka}) and (\ref{E:star}) the elements of  $X_\psi$ are $L_{\pi(A)}$-structures. Moreover the  functional correspondence $X_{\phi^\pi}\ni {\cal M}\mapsto {\cal M}^\pi\in X_\psi$ is 1-1, so $X_\psi$ is a proper class  since  $X_{\phi^\pi}$ is so. Since   VP is true  in $V$, we have
\begin{equation} \label{E:assum4}
V\models (\exists x,y)(\psi(x) \wedge \psi(y) \wedge x\neq y \wedge x\precsim_{L_{\pi(A)}} y).
\end{equation}
Pick two distinct  structures  ${\cal M}^\pi$, ${\cal N}^\pi$ of $X_\psi$ such that $${\cal M}^\pi\precsim_{L_{\pi(A)}} {\cal N}^\pi.$$
Then, by lemma \ref{L:key2}, it follows that
$V_\pi\models{\cal M}\precsim_{L_A} {\cal N}$ and also ${\cal M}\neq {\cal N}\in X_{\phi^\pi}$.
Therefore
$$V_\pi\models {\cal M}\neq {\cal N}\in X_{\phi} \wedge {\cal M}\precsim_{L_A} {\cal N}.$$
But this is the required conclusion  (\ref{E:concl}). The proof is complete. \telos

\section{VP and Choice}
What still remains open with respect to VP and the axioms of ZFC, is the  relationship of VP   with AC, namely the following:

\begin{Quest} \label{Q:ask1}
Assuming that ${\rm ZFC+VP}$ is consistent, is {\rm AC} independent from ${\rm ZF+VP}$?
\end{Quest}
Below we make two comments, one concerning  the independence of AC and one concerning   the opposite direction.

\subsection{VP and symmetric models}

We guess that AC is independent from ZF+VP. To  establish this, however,  the  most natural way  seems to be  through  the technique of permutation models of ZFA or  symmetric generic models of ZF, which are the standard tools for refuting AC.   The  technical details of  the  method  can be found in \cite{Je03}.  Also for a comprehensive list  of various symmetric models and their applications  one can  consult \cite[Part III]{HR98}.

These methods  lead inevitably to the following steps. We start with  a model $V$ of ZFA+AC+VP (or ZFC+VP) and choose a  symmetric model $HS\subset V$ for which we intend to show that   $HS\models {\rm VP}+\neg {\rm AC}$. Assuming that already $HS\models \neg {\rm AC}$,  it remains to establish that $HS\models {\rm VP}$. Let $X_\phi$ be a proper class of $L_A$-structures,  in the sense of $HS$,  for some $A\in HS$. We have to show that
\begin{equation} \label{E:req}
HS\models (\exists x,y)(\phi(x) \wedge \phi(y) \wedge x\neq y \wedge x\precsim y).
\end{equation}
Since $HS$ is an inner submodel of $V$, $X_{\phi^{HS}}$ is a proper class of $L_A$-structures in the sense of $V$, where $\phi^{HS}$ is the usual relativization of $\phi$ to $HS$. Since VP holds in $V$, we have
\begin{equation} \label{E:trans}
V\models (\exists x,y)(\phi^{HS}(x) \wedge \phi^{HS}(y) \wedge x\neq y \wedge x\precsim y).
\end{equation}
Thus the proof of the independence of AC amounts to showing that (\ref{E:req}) {\em can be  derived from} (\ref{E:trans}).

In fact proving the derivation  (\ref{E:trans})$\rightarrow$(\ref{E:req}) is a challenging problem that  cannot be settled  in the ``easy way''. The easy way would be the deduction of (\ref{E:trans})$\rightarrow$(\ref{E:req}) through  an  implication of the following form: For any language $L_A$ and any $L_A$-structures ${\cal M}$, ${\cal N}$ in $HS$,
\begin{equation} \label{E:mistakenly}
V\models {\cal M}\precsim {\cal N} \ \rightarrow \   HS\models {\cal M}\precsim {\cal N}.
\end{equation}
Obviously if (\ref{E:mistakenly}) were true for every $L_A$, the implication (\ref{E:trans})$\rightarrow$ (\ref{E:req}) would be true as well. But (\ref{E:mistakenly}) is false in general. For if we take  $L_A=\emptyset$, the $L_A$-structures are just sets and elementary embeddings are simple injections.  So  (\ref{E:mistakenly}) would imply in particular, that for all  $x,y\in HS$
\begin{equation} \label{E:mistakenly1}
V\models |x|\leq  |y| \ \rightarrow \   HS\models |x|\leq  |y|.
\end{equation}
 But since $V\models {\rm AC}$,
$$V\models (\forall x,y)(|x|\leq  |y|\vee |y|\leq  |x|).$$ So by  (\ref{E:mistakenly1}),
$$HS\models (\forall x,y)(|x|\leq  |y|\vee |y|\leq  |x|).$$
The last sentence  says that in $HS$ the cardinalities of all sets are comparable, and this is well-known to be  equivalent to AC (see for example \cite[Theorem 3.1]{RR85}). Therefore  $HS\models {\rm AC}$ which is false!

Summing up: Answering Question \ref{Q:ask1} in the affirmative amounts to finding a symmetric model $HS$ and a non-straightforward  proof of the implication (\ref{E:trans})$\rightarrow$(\ref{E:req}), for any $L_A\in HS$ and any proper class $X_\phi$ of $L_A$-structures.

\subsection{Weak forms of choice related to VP}

Now let us have a look at the opposite direction. Despite the fact  that ZF+VP is unlikely to prove AC, weaker forms of AC might be derived.

In the previous Comment  we mentioned the well-known equivalence of AC with the fact that the cardinalities of any two sets are comparable. This last formulation of AC admits natural weakenings, and one extreme such weakening is a consequence of VP.

To facilitate discussion let us say that, in ZF,  two sets $x,y$ are {\em comparable} if their cardinalities are so, that is, if either $|x|\leq |y|$ or $|y|\leq |x|$. Otherwise they are said to be {\em incomparable.}  For each formula $\phi(x)$ in one free variable, consider the following {\em comparability} axiom:
$$(Comp_\phi)  \quad \mbox{\em If $X_\phi$ is a proper class of sets, then it contains }$$
$$\mbox{\em at least two comparable elements}.$$
Then the following fact shows a slight dependence of VP with Choice.

\begin{Fac} \label{F:slight}
For every property $\phi(x)$, ${\rm ZF+VP}\vdash  Comp_\phi$.
\end{Fac}

{\em Proof.} Let $X_\phi$ be a proper class of sets. As already mentioned in  the argument in Comment 1 above, $X_\phi$ can be thought as  a proper class of $L$-structures for $L=\emptyset$, for which 1-1 mappings are elementary embeddings. Thus by VP, there are $x\neq y\in X_\phi$ such that $|x|\leq |y|$. \telos

\vskip 0.2in

Once  $Comp_\phi$ are defined, other similar axioms come up naturally. For  every standard cardinal number\footnote{By standard cardinal number we mean an initial ordinal.} $\kappa\geq 2$, consider the following comparability axiom:
$$(Comp_\kappa) \quad \quad (\forall x)[|x|=\kappa \rightarrow (\exists y\neq z\in x)(\mbox{$y$, $z$ are comparable})].$$
[$Comp_\kappa$ says that every set of cardinality $\kappa$ contains distinct comparable elements.]   $Comp_\kappa$ becomes weaker and weaker as $\kappa$ increases. Moreover $Comp_\phi$ look like   ``weakest limits'' of $Comp_\kappa$, although, apart from the implications  $Comp_n \rightarrow Comp_\phi$, $n\in\omega$, that obviously hold in  ZF, it was unknown for which $\kappa\geq \omega$ (if any) $Comp_\kappa\rightarrow Comp_\phi$ are also true.

However quite recently Lefteris Tachtsis \cite{Tach14} proved  that $Comp_\omega$ implies (over ZF)  that  finite sets coincide with Dedekind-finite sets (where $X$ is Dedekind-finite if $\aleph_0\not\leq |X|$ or,  equivalently, if there is no injection $f:X\rightarrow X$ such that $f``X\subsetneq X$). Let \textsf{F} and \textsf{DF} denote the classes of finite and Dedekind-finite sets, respectively. It is well-known that over ZF,  $\textsf{F}\subseteq \textsf{DF}$, while AC implies $\textsf{F}=\textsf{DF}$.

\begin{Thm} \label{T:Tachtsis}
{\rm (Tachtsis \cite{Tach14})} Over ${\rm ZF}$, $Comp_\omega$ implies {\rm \textsf{F}=\textsf{DF}}.
\end{Thm}

By \ref{T:Tachtsis} one can  show the following.

\begin{Prop} \label{P:omegaphi}
In {\rm ZF}, the following holds: for every $\phi$, if  $X_\phi$ is a proper class, then $Comp_\omega\rightarrow Comp_\phi$.
\end{Prop}

{\em Proof.} Assume $Comp_\omega$ is true and that $X_\phi$ is a proper class. Let us write $X$ instead of $X_\phi$. In view of $Comp_\omega$, to show that $Comp_\phi$ is true it suffices to show that there is a  set $x\subset X$ such that $|x|=\aleph_0$. Define the function $f:\omega\rightarrow On$ inductively as follows: $f(0)=0$, $f(n+1)=\min\{\beta>f(n):X\cap V_{f(n)}\subsetneq X\cap V_\beta\}$. Since $X$ is a proper class, $f$ is defined for every $n\in \omega$. (Otherwise, there exists $k\in\omega$ such that for every $\beta>f(k)$, $X\cap V_{f(k)}=X\cap V_\beta$. But then $X=X\cap V=X\cap V_{f(k)}$, so $X$ is a set.) Let $y_n=X\cap V_{f(n)}$. Each $y_n$ is a set and $y_n\subsetneq y_{n+1}$. So clearly $\bigcup_ny_n$ is an infinite set. By Theorem \ref{T:Tachtsis}, $\bigcup_ny_n$ is Dedekind-infinite, thus there is $x\subseteq \bigcup_ny_n\subset X$ such that $|x|=\aleph_0$. \telos

\vskip 0.2in

In fact the axioms  $Comp_\kappa$ are not entirely new.  $Comp_2$ is the already mentioned equivalent of AC,  that any two sets are comparable.  But also the axioms $Comp_n$, for $2\leq n<\omega$, are considered in   \cite[p. 22]{RR85}, under the name  T3($n$), and are attributed  to A. Tarski (1964). Moreover  it is shown  in  \cite[Theorem 3.4]{RR85} the significant result that $Comp_n$ is equivalent to AC, for every $2\leq n<\omega$.\footnote{I am indebted  to Lefteris Tachtsis for bringing this reference to my attention.}

The consistency of $\neg Comp_\kappa$, for every $\kappa\geq 2$, and $\neg Comp_\phi$, for {\em some} $\phi$, can be shown by the  following result of \cite{Je73}:

\begin{Thm} \label{T:Jepos}
{\rm (\cite[Thm. 11.1]{Je73})} Let $V$ be a model of  ${\rm ZFA+AC}$, with a set $A$ of atoms,  and let $\langle I,\preceq\rangle$ be a partially ordered set in $V$ such that $|A|=|I|\cdot\aleph_0$. Then there is a permutation  model $HS\subset V$   satisfying the following:
$$\mbox{There exists a family of sets $\{S_i:i\in I\}$ such that for all $i,j\in I$,}$$
\begin{equation} \label{E:pos}
i\preceq j \leftrightarrow |S_i|\leq |S_j|.
\end{equation}
\end{Thm}
The method of permutation models works also if the atoms form   a proper class rather than a set (see \cite[p. 139]{Je73}). Starting with such a model it is not difficult to strengthen \ref{T:Jepos} as follows:

\begin{Thm} \label{T:Jeposclass}
Let $V$ be a model of  ${\rm ZFA+AC}$, where now $A$ is proper class of atoms.  Let also $\langle I,\preceq\rangle$ be  a partially ordered proper class such that for every $i\in I$, $\{j:j\preceq i\}$ is a set, and $A=\{a_{in}:i\in I, n\in\omega\}$. Then there is a permutation model $HS\subset V$ and a proper class $\{S_i:i\in I\}\subset HS$ such that (\ref{E:pos}) holds.
\end{Thm}

{\em Proof.} (Sketch) By assumption for every $i\in I$, $\hat{i}=\{j\in I:j\preceq i\}$ is a set and $i\preceq j\leftrightarrow \hat{i}\subseteq \hat{j}$. Thus $\langle I,\preceq\rangle $ is embedded in $\langle{\cal P}(I),\subseteq\rangle$, where ${\cal P}(I)$ is the class of subsets of $I$, so it suffices to show that there is a class $\{S_x:x\subset I\}$ such that $x\subseteq y\leftrightarrow |S_x|\leq |S_y|$, for all $x,y\subset I$. By assumption $A=\{a_{in}:i\in I,n\in\omega\}$. For each  $x\in {\cal P}(I)$, let $S_x=\{a_{in}:i\in x,n\in\omega\}$. Consider the class-group ${\cal G}$ of permutations $\pi$ of $A$ such that $\pi(S_{\{i\}})=S_{\{i\}}$, that is, for every $i\in I$ and every $n\in\omega$ there is $m\in\omega$ such that $\pi(a_{in})=a_{im}$ (${\cal G}$ is a class of coded classes). Also let ${\cal F}$ be the filter of (suitably coded) subgroups of ${\cal G}$ generated by the ideal of finite subsets of $A$. Then the proof goes as for  Theorem \ref{T:Jepos} above.  \telos

\vskip 0.2in

Now by taking $\langle I,\preceq\rangle$ to be an antichain, that is, $i\preceq i\leftrightarrow i=j$, we obtain as a corollary to the preceding theorems the relative consistency of $\neg Comp_\kappa$, and $\neg Comp_\phi$ for a specific $\phi$.

\begin{Cor} \label{C:antichain}
(i) Let $V\models{\rm ZFA+AC}$, with a set $A$ of atoms, and let $\kappa$ be a cardinal number such that $|A|=\kappa\cdot \aleph_0$. Then there is a permutation model $HS\subset V$ such that $HS\models \neg Comp_\kappa$.

(ii) Let $V\models{\rm ZFA+AC}$, with a  proper class of atoms $A=\{a_{in}:i\in I, n\in\omega\}$. Then there is a permutation model $HS\subset V$ and a class $X_\phi=\{S_i:i\in I\}$ in $HS$ such that $HS\models\neg Comp_\phi$.
\end{Cor}

{\em Proof.} (i) Taking $I$ to be an antichain  such that $|I|=\kappa$, the claim  follows immediately from  Theorem \ref{T:Jepos}.

(ii) Similarly this  follows from Theorem \ref{T:Jeposclass} for an antichain   $I$ which is a proper class and   $X_\phi=\{S_i:i\in I\}$. \telos

\vskip 0.2in

We close with two questions concerning the relative strength of comparability axioms and their relationship with VP.

\begin{Quests} \label{Q:new}
(1) Is any of the  implications $Comp_\kappa\rightarrow Comp_\lambda$, for $\omega\leq \kappa<\lambda$,  reversible over {\rm ZF}, or over ${\rm ZF+VP}$?

(2) Does ${\rm ZF+VP}$ prove  $Comp_\kappa$, for some $\kappa\geq 2$?
\end{Quests}

{\bf Acknowledgement.} Many thanks to the anonymous referee for his very careful checking of the manuscript and the many corrections and valuable suggestions.

\end{document}